\documentclass[12pt]{article}
\usepackage{amssymb,amsmath,amsthm,amsfonts,enumerate, bbm, mathdots}

\usepackage{latexsym}
\usepackage{amscd}
\usepackage{stmaryrd}
\usepackage{color}
\usepackage[all]{xypic}
\usepackage{epsfig}
\usepackage{graphics}
\usepackage{ifthen}
\usepackage{varioref}
\usepackage{rotating}
\usepackage{extarrows}
\usepackage{cite}

\numberwithin{equation}{section}

\theoremstyle{plain}
\newtheorem{thm}{Theorem}[section]
\newtheorem{cor}[thm]{Corollary}
\newtheorem{lem}[thm]{Lemma}
\newtheorem{prop}[thm]{Proposition}
\newtheorem{rem}[thm]{Remark}

\definecolor{darkgreen}{rgb}{0.0625,0.64,0.0625}
\usepackage[
            pdfstartview=FitH,
            CJKbookmarks=true,
            bookmarksnumbered=true,
            bookmarksopen=true,
            colorlinks,
            pdfborder=001,
            linkcolor=blue,
            anchorcolor=green,
            citecolor=red
            ]{hyperref}

\newfont{\scyr}{wncyr10 scaled 550}

\def\wt{\operatorname{wt}}
\def\dep{\operatorname{dep}}
\def\height{\operatorname{-ht}}

\def\proof{\noindent {\bf Proof.\;}}

\allowdisplaybreaks

\begin{document}

\title{Sum of interpolated finite multiple harmonic $q$-series}

\date{\small ~ \qquad\qquad School of Mathematical Sciences, Tongji University \newline No. 1239 Siping Road,
Shanghai 200092, China}

\author{Zhonghua Li \thanks{E-mail address: zhonghua\_li@tongji.edu.cn}~\thanks{The authors thank the anonymous referee for his/her constructive comments. The first author is supported by the National Natural Science Foundation of
China (Grant No. 11471245).} \quad and \quad Ende Pan\thanks{E-mail address: 13162658945@163.com}}

\maketitle

\begin{abstract}
We define and study the interpolated finite multiple harmonic $q$-series. A generating function of the sums of the interpolated finite multiple harmonic $q$-series with fixed weight, depth and $i$-height is computed. Some Ohno-Zagier type relation with corollaries and some evaluation formulas of the interpolated finite multiple harmonic $q$-series at roots of unity are given.
\end{abstract}

{\small
{\bf Keywords} interpolated finite multiple harmonic $q$-series, interpolated finite $q$-multiple polylogarithms, difference equation

{\bf 2010 Mathematics Subject Classification} 11M32; 39A13
}


\section{Introduction}\label{Sec:Intro}

To study the connection between the finite and the symmetric multiple zeta values, H. Bachmann, Y. Takeyama and K. Tasaka dealt with the special values of the finite multiple harmonic $q$-series at roots of unity in \cite{BTT17}. Let $n\in \mathbb{N}$ be fixed, where $\mathbb{N}$ is the set of positive integers. Let $q$ be a complex number satisfying $q^m\neq 1$ for any $m\in\mathbb{N}$ with $m<n$. For any multi-index $\mathbf{k}=(k_1,\ldots,k_l)\in\mathbb{N}^l$ with $l\in\mathbb{N}$, the finite multiple harmonic $q$-series and their star-versions are defined by
$$z_n(\mathbf{k};q)=z_n(k_1,\ldots,k_l;q)=\sum\limits_{n>m_1>\cdots>m_l>0}\frac{q^{(k_1-1)m_1+\cdots+(k_l-1)m_l}}{[m_1]^{k_1}\cdots[m_l]^{k_l}}$$
and
$$z_n^{\star}(\mathbf{k};q)=z_n^{\star}(k_1,\ldots,k_l;q)=\sum\limits_{n>m_1\geqslant \cdots\geqslant m_l>0}\frac{q^{(k_1-1)m_1+\cdots+(k_l-1)m_l}}{[m_1]^{k_1}\cdots[m_l]^{k_l}},$$
respectively. Here for any $m\in\mathbb{N}$, $[m]$ is the $q$-integer $[m]=\frac{1-q^m}{1-q}$.  It was proved in \cite[Theorems 1.1, 1.2]{BTT17} that the values $z_n(\mathbf{k};\zeta_n)$ and $z_n^{\star}(\mathbf{k};\zeta_n)$ closely relate to the finite and the symmetric multiple zeta (star) values, where $\zeta_n$ is a fixed primitive $n$-th root of unity. Motivated by these observations, the authors of \cite{BTT17} gave a re-interpretation of the Kaneko-Zagier conjecture, which claims that the finite multiple zeta values satisfy the same $\mathbb{Q}$-linear relations as the symmetric multiple zeta values. Here $\mathbb{Q}$ is the field of rational numbers. In a latter paper \cite{BTT18}, H. Bachmann, Y. Takeyama and K. Tasaka gave some explicit evaluations and Ohno-Zagier type relations for these series at roots of unity. Note that the finite harmonic $q$-series defined above were first studied by D. M. Bradley in \cite{Bradley}.

To study the finite multiple harmonic $q$-series and their star-versions simultaneously, we introduce the interpolated finite multiple harmonic $q$-series similarly as S. Yamamoto in \cite{Yamamoto} for multiple zeta values and N. Wakabayashi in \cite{Wakabayashi} for multiple $q$-zeta values. For a multi-index $\mathbf{k}=(k_1,\ldots,k_l)\in\mathbb{N}^l$, we define its weight and depth respectively by
$$\wt(\mathbf{k})=k_1+\cdots+k_l,\quad \dep(\mathbf{k})=l.$$
Let $t$ be a formal parameter. Then we define the interpolated finite multiple harmonic $q$-series $z_n^t(\mathbf{k};q)$ as
$$z_n^t(\mathbf{k};q)=z_n^t(k_1,\ldots,k_l;q)=\sum\limits_{\mathbf{p}}(1-q)^{k-\wt(\mathbf{p})}z_n(\mathbf{p};q)t^{l-\dep(\mathbf{p})},$$
where $k=\wt(\mathbf{k})$ and $\sum\limits_{\mathbf{p}}$ is the sum where $\mathbf{p}$ runs over all multi-indices of the form $\mathbf{p}=(k_1\Box\cdots\Box k_l)$, in which each $\Box$ is filled by \lq\lq ,", \lq\lq $+$" or \lq\lq $-1+$". It is easy to see that
$$z_n^0(\mathbf{k};q)=z_n(\mathbf{k};q),\quad z_n^1(\mathbf{k};q)=z_n^{\star}(\mathbf{k};q).$$
Similar as in \cite{BTT18}, it is more convenient to study the modified versions of these $q$-series. We define
\begin{align*}
&\overline{z}_n(\mathbf{k};q)=\overline{z}_n(k_1,\ldots,k_l;q)=\sum\limits_{n>m_1>\cdots>m_l>0}\frac{q^{(k_1-1)m_1+\cdots+(k_l-1)m_l}}{(1-q^{m_1})^{k_1}\cdots(1-q^{m_l})^{k_l}},\\
&\overline{z}_n^{\star}(\mathbf{k};q)=\overline{z}_n^{\star}(k_1,\ldots,k_l;q)=\sum\limits_{n>m_1\geqslant \cdots\geqslant m_l>0}\frac{q^{(k_1-1)m_1+\cdots+(k_l-1)m_l}}{(1-q^{m_1})^{k_1}\cdots(1-q^{m_l})^{k_l}},
\end{align*}
and their interpolation polynomial
$$\overline{z}_n^t(\mathbf{k};q)=\overline{z}_n^t(k_1,\ldots,k_l;q)=\sum\limits_{\mathbf{p}}\overline{z}_n(\mathbf{p};q)t^{l-\dep(\mathbf{p})},$$
where $\mathbf{p}$ runs over all multi-indices of the form $\mathbf{p}=(k_1\Box\cdots\Box k_l)$ in which each $\Box$ is filled by \lq\lq ,", \lq\lq $+$" or \lq\lq $-1+$" as above. Then one has
$$\overline{z}_n^0(\mathbf{k};q)=\overline{z}_n(\mathbf{k};q),\quad \overline{z}_n^1(\mathbf{k};q)=\overline{z}_n^{\star}(\mathbf{k};q).$$

In this paper, we study the sums of the interpolated finite multiple harmonic $q$-series with fixed weight, depth and $1$-height, $2$-height, $\ldots$, $r$-height, which generalizes the results of \cite{BTT18}. Here for an integer $i\in\mathbb{N}$, the $i$-height of a multi-index $\mathbf{k}=(k_1,\ldots,k_l)\in\mathbb{N}^l$ is defined by
$$i\height(\mathbf{k})=\#\{j\mid k_j\geqslant i+1\}$$
as introduced in \cite{Li08}. Note that the $1$-height of a multi-index is just the height introduced by Y. Ohno and D. Zagier in \cite{OZ}. For a positive integer $r$ and nonnegative integers $k,l,h_1,\ldots,h_r$, we set
$$G_n^t(k,l,h_1,\ldots,h_r;q)=\sum\limits_{\mathbf{k}\in I(k,l,h_1,\ldots,h_r)}\overline{z}_n^t(\mathbf{k};q),$$
where $I(k,l,h_1,\ldots,h_r)$ is the set of multi-indices of weight $k$, depth $l$, $1$-height $h_1$, $\ldots$, $r$-height $h_r$. As usual, the sum is treated as $0$ whenever the index set is empty except for $G^t_n(0,0,\ldots,0;q)=1$. For formal variables $u_1,\ldots,u_{r+2}$, we define the generating function
\begin{align*}
&\Psi_n^t(q)=\Psi_n^t(u_1,\ldots,u_{r+2};q)\\
=&\sum\limits_{k,l,h_1,\ldots,h_r\geqslant 0}G_n^t(k,l,h_1,\ldots,h_r;q)u_1^{k-l-\sum_{j=1}^rh_j}u_2^{l-h_1}u_3^{h_1-h_2}\cdots u_{r+1}^{h_{r-1}-h_r}u_{r+2}^{h_r}.
\end{align*}
Then similarly as in \cite{AOW,BTT18,Li08,Li10,LW}, we obtain the following theorem.

\begin{thm}\label{Thm:GeneratingFunction}
Let $r$ be a positive integer and $u_1,\ldots,u_{r+2}$ be variables. Let
\begin{align}
x_1=\frac{u_1}{1+u_1}
\label{Eq:x1}
\end{align}
and
\begin{align}
x_i=\sum\limits_{j=i}^{r+1}(-1)^{j-i}\binom{j-2}{i-2}\left(u_j-\frac{u_{r+2}}{u_1^{r+2-j}}\right)+\frac{u_{r+2}}{u_1^{r+2-i}(1+u_1)^{i-1}}
\label{Eq:xi}
\end{align}
for $i=2,\ldots,r+2$. We have
$$\Psi_n^t(q)=\frac{\prod\limits_{j=1}^{n-1}P^{t-1}(1-q^j)}{\prod\limits_{j=1}^{n-1}P^{t}(1-q^j)},$$
where
$$P^t(T)=T^{r+1}-(x_1+tx_2)T^r-t\sum\limits_{i=0}^{r-1}(x_{r+2-i}-x_1x_{r+1-i})T^i.$$
\end{thm}

From Theorem \ref{Thm:GeneratingFunction}, it is easy to see that
$$\Psi_n^{1-t}(u_1,-u_2,\ldots,-u_{r+2};q)=\frac{\prod\limits_{j=1}^{n-1}P^{t}(1-q^j)}{\prod\limits_{j=1}^{n-1}P^{t-1}(1-q^j)}.$$
Hence we find
$$\Psi_n^t(u_1,\ldots,u_{r+2};q)\Psi_n^{1-t}(u_1,-u_2,\ldots,-u_{r+2};q)=1.$$
Let $t=0$, we get a relation between the generating function of the sums of the finite multiple harmonic $q$-series and the generating function of their star-versions. See also \cite[Theorem 1.2]{BTT18}. And setting $t=\frac{1}{2}$, we get
$$\Psi_n^{\frac{1}{2}}(u_1,\ldots,u_{r+2};q)\Psi_n^{\frac{1}{2}}(u_1,-u_2,\ldots,-u_{r+2};q)=1,$$
which would give some relations of the sums $G_n^{\frac{1}{2}}(k,l,h_1,\ldots,h_r;q)$. We omit the formulas.

Now let $q=\zeta_n$ be a fixed primitive $n$-th root of unity. Theorem \ref{Thm:GeneratingFunction} implies the following corollary.

\begin{cor}\label{Cor:Sum-rational}
For any nonnegative integers $k,l,h_1,\ldots,h_r$, we have
$$G^t_n(k,l,h_1,\ldots,h_r;\zeta_n)\in\mathbb{Q}[t].$$
\end{cor}

In the case of  $r=1$, $G^t_n(k,l,s;q)$ is the sum of the modified interpolated finite multiple harmonic $q$-series with fixed weight, depth and height. We have the following Ohno-Zagier type relation.

\begin{thm}\label{Thm:OhnoZagier-tqFMHS}
Let $u_1,u_2,u_3$ be formal variables. Then we have
$$\Psi_n^t(u_1,u_2,u_3;\zeta_n)=\sum\limits_{k,l,s\geqslant 0}G_n^t(k,l,s;\zeta_n)u_1^{k-l-s}u_2^{l-s}u_3^s
=\frac{U_n^{t-1}(u_1,u_2,u_3)}{U_n^t(u_1,u_2,u_3)},
$$
where
\begin{align*}
U_n^t(u_1,u_2,u_3)=&\sum\limits_{a,b\geqslant 0\atop a+b\leqslant n-1}\frac{1}{n-a-b}\binom{n-a-1}{b}\binom{n-b-1}{a}\\
&\quad\times t^{n-a-b-1}(1+u_1)^a(1-tu_2)^b(u_3-u_1u_2)^{n-a-b-1}.
\end{align*}
\end{thm}

Moreover, setting $u_3=u_1u_2$ in Theorem \ref{Thm:OhnoZagier-tqFMHS}, we meet the sums of the modified interpolated finite multiple harmonic $q$-series with fixed weight and depth. For nonnegative integers $k$ and $l$, let $I(k,l)$ be the set of multi-indices of weight $k$ and depth $l$ and
$$G_n^t(k,l;q)=\sum\limits_{\mathbf{k}\in I(k,l)}\overline{z}_n^{t}(\mathbf{k};q).$$
We have two expressions for the sum $G_n^t(k,l;\zeta_n)$ as displayed in the following corollary, in which \eqref{Eq:SumFormula-2} may be regarded as the sum formula for the interpolated finite multiple harmonic $q$-series at roots of unity.

\begin{cor}\label{Cor:SumFormula-tqFHMS}
For integers $k,l$ with $k\geqslant l\geqslant 0$, we have
\begin{align}
G_n^t(k,l;\zeta_n)=\sum\limits_{m=0}^k\frac{(-1)^m}{n^{m+1}}\sum\limits_{{{i_0+\cdots+i_m=l\atop j_0+\cdots+j_m=k}\atop {0\leqslant i_a\leqslant j_a\leqslant n-1\atop a=0,\ldots,m}}\atop j_1,\ldots,j_m\geqslant 1}\prod\limits_{a=0}^m\binom{n}{j_a+1}(1-t)^{i_0}(-t)^{l-i_0}
\label{Eq:SumFormula-1}
\end{align}
and
\begin{align}
G_n^t(k,l;\zeta_n)=&-\sum\limits_{m=0}^k\frac{1}{n^{m+1}}\sum\limits_{{i_0+\cdots+i_m=l\atop l_0+\cdots+l_m+j_0+\cdots+j_m=k}\atop{{0\leqslant i_0\leqslant j_0\leqslant n-1\atop 1\leqslant i_a\leqslant j_a\leqslant n-1}\atop {a=1,\ldots,m\atop l_0,\ldots,l_m\geqslant 0}}}\prod\limits_{a=0}^m\binom{n}{j_a+1}\overline{z}_n(l_a;\zeta_n)\nonumber\\
&\qquad\times(1-t)^{i_0}(-t)^{l-i_0}.
\label{Eq:SumFormula-2}
\end{align}
Here we use the convention $\overline{z}_n(0;\zeta_n)=-1$.
\end{cor}

Setting $t=0$ in \eqref{Eq:SumFormula-2}, we have
\begin{align}
\sum\limits_{\mathbf{k}\in I(k,l)}\overline{z}_n(\mathbf{k};\zeta_n)=-\frac{1}{n}\sum\limits_{j=l}^k\binom{n}{j+1}\overline{z}_n(k-j;\zeta_n),\quad (k\geqslant l\geqslant 0),
\label{Eq:SumFormula-FMHS}
\end{align}
which is equivalent to \cite[(3.14)]{BTT18}. In fact, putting $l=1$ in \eqref{Eq:SumFormula-FMHS}, we have
$$\sum\limits_{j=0}^k\binom{n}{j+1}\overline{z}_n(k-j;\zeta_n)=0,\quad (k\geqslant 1).$$
Then for $k\geqslant l\geqslant 1$, the equation \eqref{Eq:SumFormula-FMHS} is equivalent to
$$\sum\limits_{\mathbf{k}\in I(k,l)}\overline{z}_n(\mathbf{k};\zeta_n)=\frac{1}{n}\sum\limits_{j=0}^{l-1}\binom{n}{j+1}\overline{z}_n(k-j;\zeta_n),$$
which is just \cite[(3.14)]{BTT18}.

Let $u_1=u_3=0$ in Theorem \ref{Thm:OhnoZagier-tqFMHS}, we obtain the evaluation formula of $\overline{z}_n^t(\{1\}^l;\zeta_n)$. Here and below, $\{k\}^l$ stands for $\underbrace{k,\ldots,k}_{l\text{\;terms}}$. Moreover, let $u_1=\cdots=u_{r+1}=0$ and $q=\zeta_n$ in Theorem \ref{Thm:GeneratingFunction}, we get the evaluation formulas of $\overline{z}_n^t(\{k\}^l;\zeta_n)$ for any integer $k\geqslant 2$. More explicitly, we have the following formulas.

\begin{cor}\label{Cor:Eva-k-l}
For any nonnegative integer $l$, we have
\begin{align}
\overline{z}_n^t(\{1\}^l;\zeta_n)=&\sum\limits_{m=0}^l\frac{(-1)^m}{n^{m+1}}\sum\limits_{i_0+\cdots+i_m=l\atop {0\leqslant i_0\leqslant n-1\atop 1\leqslant i_1,\ldots,i_m\leqslant n-1}}\prod\limits_{a=0}^m\binom{n}{i_a+1}(1-t)^{i_0}(-t)^{l-i_0},
\label{Eq:Eva-1-l}\\
\overline{z}_n^t(\{2\}^l;\zeta_n)=&\sum\limits_{m=0}^l\frac{(-1)^m}{n^{m+1}}\sum\limits_{i_0+\cdots+i_m=l\atop {0\leqslant i_0\leqslant n-1\atop 1\leqslant i_1,\ldots,i_m\leqslant n-1}}\prod\limits_{a=0}^m\frac{1}{i_a+1}\binom{n+i_a}{2i_a+1}(t-1)^{i_0}t^{l-i_0},
\label{Eq:Eva-2-l}\\
\overline{z}_n^t(\{3\}^l;\zeta_n)=&\sum\limits_{m=0}^l\frac{(-1)^m}{n^{2m+2}}\sum\limits_{i_0+\cdots+i_m=l\atop {0\leqslant i_0\leqslant n-1\atop 1\leqslant i_1,\ldots,i_m\leqslant n-1}}\prod\limits_{a=0}^m\frac{1}{i_a+1}\left[\binom{n+i_a}{3i_a+2}\right.\nonumber\\
&\qquad\qquad\qquad\qquad\left.+(-1)^{i_a}\binom{n+2i_a+1}{3i_a+2}\right](t-1)^{i_0}t^{l-i_0}.
\label{Eq:Eva-3-l}
\end{align}
\end{cor}

As in \cite{BTT17,BTT18}, we can consider the limit of the interpolated finite multiple $q$-series when $q=\zeta_n=e^{2\pi \sqrt{-1}/n}$ as $n\rightarrow \infty$. For any multi-index $\mathbf{k}=(k_1,\ldots,k_l)\in\mathbb{N}^l$, we set
$$\xi(\mathbf{k})=\lim\limits_{n\rightarrow\infty}z_n(\mathbf{k};e^{2\pi \sqrt{-1}/n}),$$
which is well defined by \cite[Theorem 1.2]{BTT17}. Now one can define
$$\xi^t(\mathbf{k})=\xi^t(k_1,\ldots,k_l)=\sum\limits_{\mathbf{p}}\xi(\mathbf{p})t^{l-\dep(\mathbf{p})},$$
where $\mathbf{p}$ runs over all multi-indices of the form $\mathbf{p}=(k_1\Box\cdots\Box k_l)$ in which each $\Box$ is filled by \lq\lq ," or \lq\lq $+$". It is easy to see that
$$\xi^t(\mathbf{k})=\lim\limits_{n\rightarrow \infty}z_n^t(\mathbf{k};e^{2\pi \sqrt{-1}/n}).$$
From \eqref{Eq:SumFormula-2} and Corollary \ref{Cor:Eva-k-l}, we can obtain the sum formula of $\xi^t(\mathbf{k})$ and the evaluation formulas of $\xi^t(\{k\}^l)$ for $k=1,2,3$. For example, we have
$$\xi^t(\{1\}^l)=\sum\limits_{m=0}^l(-1)^m\sum\limits_{i_0+\cdots+i_m=l\atop i_0\geqslant 0,i_1,\ldots,i_m\geqslant 1}\frac{(1-t)^{i_0}(-t)^{l-i_0}}{(i_0+1)!\cdots(i_m+1)!}(-2\pi \sqrt{-1})^l.$$
We omit other formulas and leave them to the interested readers. And one can use \cite[Theorem 1.1]{BTT17} and the results mentioned above to deduce some relations of the interpolated finite multiple zeta values defined by S. Seki in \cite{Seki}. We also omit the details.

The paper is organized as follows. In Section \ref{Sec:Sum-tqFMP}, we define the interpolated finite $q$-multiple polylogarithms, and consider similar sums of these functions. Applying the results of Section \ref{Sec:Sum-tqFMP}, we prove Theorem \ref{Thm:GeneratingFunction} in Section \ref{Sec:Sum-tqFMHS}. In Section \ref{Sec:RootUnity}, we specify $q=\zeta_n$ and prove Corollary \ref{Cor:Sum-rational}, Theorem \ref{Thm:OhnoZagier-tqFMHS}, Corollary \ref{Cor:SumFormula-tqFHMS} and Corollary \ref{Cor:Eva-k-l}.


\section{Sum of interpolated finite $q$-multiple polylogarithms} \label{Sec:Sum-tqFMP}

The proof of Theorem \ref{Thm:GeneratingFunction} is similar to that of other Ohno-Zagier type relations \cite{AOW,BTT18,Li08,Li10,LW}. We first study a similar sum of the interpolated finite $q$-multiple polylogarithms in this section. Then we get Theorem \ref{Thm:GeneratingFunction} in the next section.

For a positive integer $n$ and a multi-index $\mathbf{k}=(k_1,\ldots,k_l)\in\mathbb{N}^l$ with $l\in\mathbb{N}$, as in \cite{BTT18}, we define two polynomials in $z$ by
\begin{align*}
&L_{n,\mathbf{k}}(z;q)=\sum\limits_{n>m_1>\cdots>m_l>0}\frac{z^{m_1}}{(1-q^{m_1})^{k_1}\cdots(1-q^{m_l})^{k_l}},\\
&L_{n,\mathbf{k}}^{\star}(z;q)=\sum\limits_{n>m_1\geqslant \cdots\geqslant m_l>0}\frac{z^{m_1}}{(1-q^{m_1})^{k_1}\cdots(1-q^{m_l})^{k_l}}.
\end{align*}
Note that the degree of each of the polynomials $L_{n,\mathbf{k}}(z;q)$ and $L_{n,\mathbf{k}}^{\star}(z;q)$ is less than $n$. Then we define the interpolated polynomial of these two functions by
$$L^t_{n,\mathbf{k}}(z;q)=\sum\limits_{\mathbf{p}}L_{n,\mathbf{p}}(z;q)t^{l-\dep(\mathbf{p})},$$
where $\mathbf{p}$ runs over all multi-indices of the form $\mathbf{p}=(k_1\Box\cdots\Box k_l)$ in which each $\Box$ is filled by \lq\lq ," or \lq\lq $+$". It is easy to see
$$L^0_{n,\mathbf{k}}(z;q)=L_{n,\mathbf{k}}(z;q),\quad L^1_{n,\mathbf{k}}(z;q)=L^\star_{n,\mathbf{k}}(z;q),$$
as expected.

We define the $q$-difference operator $\Theta_q$ by
$$(\Theta_qf)(z)=f(z)-f(qz).$$
Then for $\mathbf{k}=(k_1,\ldots,k_l)\in\mathbb{N}^l$, we have
$$\Theta_qL_{n,\mathbf{k}}(z;q)=\begin{cases}
L_{n,(k_1-1,k_2,\ldots,k_l)}(z;q) & \text{if\;} k_1\geqslant 2,\\
\frac{z}{1-z}L_{n,(k_2,\ldots,k_l)}(z;q)-\frac{z^{n}}{1-z}L_{n,(k_2,\ldots,k_l)}(1;q) & \text{if\;} k_1=1, l\geqslant 2,\\
\frac{z-z^{n}}{1-z} & \text{if\;} k_1=l=1.
\end{cases}$$
From above formulas we get
\begin{align}
\Theta_qL_{n,\mathbf{k}}^t(z;q)=\begin{cases}
L_{n,(k_1-1,k_2,\ldots,k_l)}^t(z;q) & \text{if\;} k_1\geqslant 2,\\
\left(t+\frac{z}{1-z}\right)L_{n,(k_2,\ldots,k_l)}^t(z;q)&\\
\qquad\qquad -\frac{z^{n}}{1-z}L_{n,(k_2,\ldots,k_l)}^t(1;q) & \text{if\;} k_1=1, l\geqslant 2,\\
\frac{z-z^{n}}{1-z} & \text{if\;} k_1=l=1,
\end{cases}
\label{Eq:DiffEq-tqFMPL}
\end{align}
which can be proved similarly as that of \cite[Lemma 5]{LW}.

Let $r$ be a positive integer and $k,l,h_1,\ldots,h_r$ be nonnegative integers. For an integer $j$ with $-1\leqslant j\leqslant r-1$, let $I_j(k,l,h_1,\ldots,h_r)$ be the set of multi-indices $\mathbf{k}=(k_1,\ldots,k_l)$ satisfying
$$\wt(\mathbf{k})=k,\quad \dep(\mathbf{k})=l,\quad i\height(\mathbf{k})=h_i, i=1,\ldots,r$$
and $k_1\geqslant j+2$. Note that the index set $I(k,l,h_1,\ldots,h_r)$ defined in Section \ref{Sec:Intro} is just $I_{-1}(k,l,h_1,\ldots,h_r)$. Now we define sums of the interpolated finite $q$-multiple polylogarithms by
$$X_{n,j}^t(k,l,h_1,\ldots,h_r;z;q)=\sum\limits_{\mathbf{k}\in I_j(k,l,h_1,\ldots,h_r)}L_{n,\mathbf{k}}^t(z;q),$$
which is $0$ whenever the index set is empty except for $X_{n,-1}^t(0,0,\ldots,0;z;q)=1$. We also abbreviate $X_{n,-1}^t(k,l,h_1,\ldots,h_r;z;q)$ to $X_{n}^t(k,l,h_1,\ldots,h_r;z;q)$. Using \eqref{Eq:DiffEq-tqFMPL}, we can prove the following result similarly as that of \cite[Lemma 6]{LW}.

\begin{lem}\label{Lem:DifferEq-SumtqFMPL}
Let $k,l,h_1,\ldots,h_r$ be integers satisfying $k\geqslant l+\sum_{j=1}^rh_j$ and $l\geqslant h_1\geqslant \cdots\geqslant h_r\geqslant 0$.
\begin{enumerate}
  \item[(i)] If $h_r\geqslant 1$, we have
  \begin{align*}
  \Theta_qX_{n,r-1}^t(k,l,h_1,\ldots,h_r;z;q)=&X_{n,r-1}^t(k-1,l,h_1,\ldots,h_r;z;q)\\
  &+X_{n,r-2}^t(k-1,l,h_1,\ldots,h_{r-1},h_r-1;z;q)\\
  &-X_{n,r-1}^t(k-1,l,h_1,\ldots,h_{r-1},h_r-1;z;q).
  \end{align*}
  \item[(ii)] For $0\leqslant j\leqslant r-2$ with $h_{j+1}\geqslant 1$, we have
  \begin{align*}
  &\Theta_q\left[X_{n,j}^t(k,l,h_1,\ldots,h_r;z;q)-X_{n,j+1}^t(k,l,h_1,\ldots,h_r;z;q)\right]\\
  =&X_{n,j-1}^t(k-1,l,h_1,\ldots,h_j,h_{j+1}-1,h_{j+2},\ldots,h_r;z;q)\\
  &-X_{n,j}^t(k-1,l,h_1,\ldots,h_j,h_{j+1}-1,h_{j+2},\ldots,h_r;z;q).
  \end{align*}
  \item[(iii)] If $l\geqslant 2$, we have
  \begin{align*}
  &\Theta_q\left[X_{n}^t(k,l,h_1,\ldots,h_r;z;q)-X_{n,0}^t(k,l,h_1,\ldots,h_r;z;q)\right]\\
  =&\left(t+\frac{z}{1-z}\right)X_{n}^t(k-1,l-1,h_1,\ldots,h_r;z;q)\\
  &-\frac{z^n}{1-z}X_{n}^t(k-1,l-1,h_1,\ldots,h_r;1;q).
  \end{align*}
\end{enumerate}
\end{lem}

Let $x_1,\ldots,x_{r+2}$ be variables. For $-1\leqslant j\leqslant r-1$, we define the generating functions
\begin{align*}
&\Phi_j^t(z)=\Phi_j^t(z;q)=\Phi_{n,j}^t(x_1,\ldots,x_{r+2};z;q)\\
=&\sum\limits_{k,l,h_1,\ldots,h_r\geqslant 0}X_{n,j}^t(k,l,h_1,\ldots,h_r;z;q)x_1^{k-l-\sum_{i=1}^rh_i}x_2^{l-h_1}x_3^{h_1-h_2}\cdots x_{r+1}^{h_{r-1}-h_r}x_{r+2}^{h_r}.
\end{align*}
Let us denote $\Phi_{-1}^t(z)$ by $\Phi^t(z)$. Then from Lemma \ref{Lem:DifferEq-SumtqFMPL}, we get the following proposition, which can be proved similarly as that of \cite[Proposition 7]{LW}.

\begin{prop}\label{Prop:DiffEq-GeneFun-tqFMPL}
We have
$$\begin{cases}
\Theta_q\Phi_{r-1}^t(z)=x_1\Phi_{r-1}^t(z)+\frac{x_{r+2}}{x_{r+1}}\left(\Phi_{r-2}^t(z)-\Phi_{r-1}^t(z)-\delta_{r,1}\right),&\\
\Theta_q\left(\Phi_j^t(z)-\Phi_{j+1}^t(z)\right)=\frac{x_{j+3}}{x_{j+2}}\left(\Phi_{j-1}^t(z)-\Phi_{j}^t(z)\right), & j=1,\ldots,r-2,\\
\Theta_q\left(\Phi_0^t(z)-\Phi_1^t(z)\right)=\frac{x_3}{x_2}\left(\Phi^t(z)-\Phi_0^t(z)-1\right),&\\
\Theta_q\left(\Phi^t(z)-\Phi_0^t(z)\right)=\left(t+\frac{z}{1-z}\right)x_2\Phi^t(z)-tx_2-\frac{z^n}{1-z}x_2\Phi^t(1),
\end{cases}$$
where $\delta$ stands for the Kronecker's delta.
\end{prop}

Similarly as that of \cite[Corollary 8]{LW}, eliminating other generating functions from Proposition \ref{Prop:DiffEq-GeneFun-tqFMPL}, we get the $q$-difference equation satisfying by $\Phi_{r-1}^t(z)$.

\begin{cor}\label{Cor:DiffEq-Phi-r-1}
The function $y_0=\Phi_{r-1}^t(z)$ satisfies the following $q$-difference equation
$$\left(P^t(\Theta_q)-zP^{t-1}(\Theta_q)\right)y_0=zx_{r+2}-z^nx_{r+2}\Phi^t(1),$$
where
$$P^t(T)=T^{r+1}-(x_1+tx_2)T^r-t\sum\limits_{i=0}^{r-1}(x_{r+2-i}-x_1x_{r+1-i})T^i,$$
as defined in Section \ref{Sec:Intro}.
\end{cor}

Finally, we get $\Phi^t(1)$ as in the following theorem.

\begin{thm}\label{Thm:Phi-1-tqFMPL}
We have
$$\Phi^t(1)=\Phi^t(1;q)=\frac{\prod\limits_{j=1}^{n-1}P^{t-1}(1-q^j)}{\prod\limits_{j=1}^{n-1}P^t(1-q^j)}.$$
\end{thm}

\proof
Set $y_0=\Phi_{r-1}^t(z)=\sum_{i=1}^{n-1}c_iz^i$. Since
$$\Theta_qz^i=(1-q^i)z^i,$$
we get
$$P^t(\Theta_q)z^i=P^t(1-q^i)z^i.$$
Hence from Corollary \ref{Cor:DiffEq-Phi-r-1}, we have
$$\sum\limits_{i=1}^{n-1}c_iP^t(1-q^i)z^i-\sum\limits_{i=1}^{n-1}c_iP^{t-1}(1-q^i)z^{i+1}=zx_{r+2}-z^nx_{r+2}\Phi^t(1).$$
Comparing the coefficients of $z^i$, we get
$$\begin{cases}
c_1P^t(1-q)=x_{r+2},&\\
c_iP^t(1-q^i)=c_{i-1}P^{t-1}(1-q^{i-1}), & i=2,\ldots,n-1,\\
c_{n-1}P^{t-1}(1-q^{n-1})=x_{r+2}\Phi^t(1).&
\end{cases}$$
Therefore we have
\begin{align}
c_i=\frac{\prod\limits_{j=1}^{i-1}P^{t-1}(1-q^j)}{\prod\limits_{j=1}^iP^t(1-q^j)}x_{r+2}
\label{Eq:c-i}
\end{align}
for $i=1,\ldots,n-1$, which can be proved by induction on $i$. Finally we get the desired result.
\qed

\begin{rem}
Assume that
$$P^t(T)=(T+\alpha_1^t)\cdots(T+\alpha_{r+1}^t),$$
and let
$$a_i^{t}=\frac{q}{1+\alpha_i^t},\quad i=1,\ldots,r+1.$$
Then we have
$$\Phi_{r-1}^t(z)=\frac{zx_{r+2}}{P^t(1-q)}\;_{r+2}\phi_{r+1}\left[\begin{matrix}
q, a_1^{t-1},\ldots,a_{r+1}^{t-1}\\
qa_1^t,\ldots,qa_{r+1}^t
\end{matrix};q;\frac{a_1^{t}\cdots a_{r+1}^{t}}{a_1^{t-1}\cdots a_{r+1}^{t-1}}z\right]_{n-2}.$$
Here we use the notion of the truncated basic hypergeometric series defined as
$$\;_{p+1}\phi_p\left[\begin{matrix}
a_1, \ldots, a_{p+1}\\
b_1, \ldots, b_p
\end{matrix};q,z\right]_n=\sum\limits_{i=0}^n\frac{(a_1,\ldots,a_{p+1};q)_i}{(q,b_1,\ldots,b_p;q)_i}z^i,$$
where $(a_1,\ldots,a_{p+1};q)_i=(a_1;q)_i\cdots(a_{p+1};q)_i$ with the $q$-shifted factorial $(a;q)_i$ defined as
$$(a;q)_i=\begin{cases}
1 & \text{if\;} i=0,\\
(1-a)(1-aq)\cdots(1-aq^{i-1}) & \text{if\;} i=1,2,\ldots.
\end{cases}$$
Similarly as in \cite{LW}, one can represent $\Phi_j^t(z)$ by the truncated basic hypergeometric series for any $j$.
\end{rem}


\section{Sum of interpolated finite multiple harmonic $q$-series}\label{Sec:Sum-tqFMHS}

Using Theorem \ref{Thm:Phi-1-tqFMPL}, one can prove Theorem \ref{Thm:GeneratingFunction}. For that purpose, we prepare the following lemma, which represents $L_{n,\mathbf{k}}^t(1;q)$ by the interpolated finite multiple harmonic $q$-series.

\begin{lem}
For any multi-index $\mathbf{k}=(k_1,\ldots,k_l)\in\mathbb{N}^l$, we have
\begin{align}
L_{n,\mathbf{k}}^t(1;q)=\sum\limits_{a_1=1}^{k_1}\cdots\sum\limits_{a_l=1}^{k_l}\binom{k_1-1}{a_1-1}\cdots\binom{k_l-1}{a_l-1}\overline{z}_n^t(a_1,\ldots,a_l;q).
\label{Eq:tqFMPL-tqFMHS}
\end{align}
\end{lem}

\proof
For any multi-index $(a_1,\ldots,a_l)\in\mathbb{N}^l$, the multi-index $\mathbf{p}$ in the definition of $\overline{z}_n^t(a_1,\ldots,a_l;q)$ has the form
$$\mathbf{p}=(a_1\Box\cdots\Box a_{i_1},a_{i_1+1}\Box\cdots\Box a_{i_1+i_2},\ldots,a_{i_1+\cdots+i_{c-1}+1}\Box\cdots\Box a_{i_1+\cdots+i_c}),$$
in which each $\Box$ is filled by \lq\lq $+$" or \lq\lq $-1+$". Let $m_j$ be the numbers of \lq\lq $-1+$" for $j=1,\ldots,c$. Then we have $0\leqslant m_j\leqslant i_j-1$, and
\begin{align*}
\overline{z}_n^t(a_1,\ldots,a_l;q)=&\sum\limits_{c=1}^l\sum\limits_{i_1+\cdots+i_c=l\atop i_1,\ldots,i_c\geqslant 1}\sum\limits_{0\leqslant m_j\leqslant i_j-1\atop j=1,\ldots,c}\binom{i_1-1}{m_1}\cdots\binom{i_c-1}{m_c}\\
&\times\overline{z}_n(a_{I_1}-m_1,\ldots,a_{I_c}-m_c;q)t^{l-c},
\end{align*}
where
$$\begin{cases}
a_{I_1}=a_1+\cdots+a_{i_1},&\\
a_{I_2}=a_{i_1+1}+\cdots+a_{i_1+i_2},&\\
\qquad\qquad\vdots\\
a_{I_c}=a_{i_1+\cdots+i_{c-1}+1}+\cdots+a_{i_1+\cdots+i_c}.&
\end{cases}$$
Hence the right-hand side of \eqref{Eq:tqFMPL-tqFMHS} is
\begin{align*}
&\sum\limits_{a_1=1}^{k_1}\cdots\sum\limits_{a_l=1}^{k_l}\left\{\prod\limits_{j=1}^l\binom{k_j-1}{a_j-1}\right\}\sum\limits_{c=1}^l\sum\limits_{i_1+\cdots+i_c=l\atop i_1,\ldots,i_c\geqslant 1}\sum\limits_{0\leqslant m_j\leqslant i_j-1\atop j=1,\ldots,c}\left\{\prod\limits_{j=1}^c\binom{i_j-1}{m_j}\right\}\\
&\quad\times\overline{z}_n(a_{I_1}-m_1,\ldots,a_{I_c}-m_c;q)t^{l-c}\\
=&\sum\limits_{c=1}^l\sum\limits_{i_1+\cdots+i_c=l\atop i_1,\ldots,i_c\geqslant 1}\sum\limits_{0\leqslant m_j\leqslant i_j-1\atop j=1,\ldots,c}\left\{\prod\limits_{j=1}^c\binom{i_j-1}{m_j}\right\}\sum\limits_{a_1=1}^{k_1}\cdots\sum\limits_{a_l=1}^{k_l}\left\{\prod\limits_{j=1}^l\binom{k_j-1}{a_j-1}\right\}\\
&\quad\times\overline{z}_n(a_{I_1}-m_1,\ldots,a_{I_c}-m_c;q)t^{l-c}.
\end{align*}
We first compute the inner sums
$$\sum\limits_{a_1=1}^{k_1}\cdots\sum\limits_{a_l=1}^{k_l}\left\{\prod\limits_{j=1}^l\binom{k_j-1}{a_j-1}\right\}\overline{z}_n(a_{I_1}-m_1,\ldots,a_{I_c}-m_c;q),$$
which is equal to
\begin{align*}
&\prod\limits_{j=1}^c\left\{\sum\limits_{b_j=i_j}^{k_{I_j}}\sum\limits_{a_{i_1+\cdots+i_{j-1}+1}+\cdots+a_{i_1+\cdots+i_j}=b_j\atop a_{\bullet}\geqslant 1}\prod\limits_{d=i_1+\cdots+i_{j-1}+1}^{i_1+\cdots+i_{j}}\binom{k_{d}-1}{a_{d}-1}\right\}\\
&\quad\times \overline{z}_n(b_1-m_1,\ldots,b_c-m_c;q)\\
=&\prod\limits_{j=1}^c\left\{\sum\limits_{b_j=i_j}^{k_{I_j}}\binom{k_{I_j}-i_j}{b_j-i_j}\right\}\overline{z}_n(b_1-m_1,\ldots,b_c-m_c;q)\\
=&\prod\limits_{j=1}^c\left\{\sum\limits_{b_j=i_j-m_j}^{k_{I_j}-m_j}\binom{k_{I_j}-i_j}{b_j+m_j-i_j}\right\}\overline{z}_n(b_1,\ldots,b_c;q),
\end{align*}
where
$$\begin{cases}
k_{I_1}=k_1+\cdots+k_{i_1},&\\
k_{I_2}=k_{i_1+1}+\cdots+k_{i_1+i_2},&\\
\qquad\qquad\vdots\\
k_{I_c}=k_{i_1+\cdots+i_{c-1}+1}+\cdots+k_{i_1+\cdots+i_c}.&
\end{cases}$$
Thus the right-hand side of \eqref{Eq:tqFMPL-tqFMHS} is
\begin{align*}
&\sum\limits_{c=1}^l\sum\limits_{i_1+\cdots+i_c=l\atop i_1,\ldots,i_c\geqslant 1}\prod\limits_{j=1}^c\left\{\sum\limits_{m_j=0}^{i_j-1}\sum\limits_{b_j=i_j-m_j}^{k_{I_j}-m_j}\binom{i_j-1}{m_j}\binom{k_{I_j}-i_j}{k_{I_j}-b_j-m_j}\right\}\overline{z}_n(\mathbf{b};q)t^{l-c}\\
=&\sum\limits_{c=1}^l\sum\limits_{i_1+\cdots+i_c=l\atop i_1,\ldots,i_c\geqslant 1}\prod\limits_{j=1}^c\left\{\sum\limits_{b_j=1}^{k_{I_j}}\sum\limits_{m_j=0}^{i_j-1}\binom{i_j-1}{m_j}\binom{k_{I_j}-i_j}{k_{I_j}-b_j-m_j}\right\}\overline{z}_n(\mathbf{b};q)t^{l-c}\\
=&\sum\limits_{c=1}^l\sum\limits_{i_1+\cdots+i_c=l\atop i_1,\ldots,i_c\geqslant 1}\prod\limits_{j=1}^c\left\{\sum\limits_{b_j=1}^{k_{I_j}}\binom{k_{I_j}-1}{b_j-1}\right\}\overline{z}_n(\mathbf{b};q)t^{l-c},
\end{align*}
where $\mathbf{b}=(b_1,\ldots,b_c)$. By \cite[(3.4)]{BTT18}, we have
$$L_{n,(k_{I_1},\ldots,k_{I_c})}(1;q)=\prod\limits_{j=1}^c\left\{\sum\limits_{b_j=1}^{k_{I_j}}\binom{k_{I_j}-1}{b_j-1}\right\}\overline{z}_n(\mathbf{b};q).$$
Therefore the right-hand side of \eqref{Eq:tqFMPL-tqFMHS} is
$$\sum\limits_{c=1}^l\sum\limits_{i_1+\cdots+i_c=l\atop i_1,\ldots,i_c\geqslant 1}L_{n,(k_{I_1},\ldots,k_{I_c})}(1;q)t^{l-c},$$
which is just $L_{n,\mathbf{k}}^t(1;q)$. We finish the proof of \eqref{Eq:tqFMPL-tqFMHS}.
\qed

\begin{rem}
The proof of \eqref{Eq:tqFMPL-tqFMHS} offered here works for \cite[Lemma 16]{LW} in the case of the interpolated multiple $q$-zeta values, where the algebraic method was  used.
\end{rem}

Using \eqref{Eq:tqFMPL-tqFMHS}, we represent $\Phi^t(1;q)$ by $\Psi^t_n(q)$ defined in Section \ref{Sec:Intro} as in the following lemma.

\begin{lem}\label{Lem:Phi-1-Psi}
For variables $x_1,\ldots,x_{r+2}$, let
\begin{align}
u_1=\frac{x_1}{1-x_1}
\label{Eq:u1}
\end{align}
and
\begin{align}
u_i=\sum\limits_{j=i}^{r+1}\binom{j-2}{i-2}\left(x_j-\frac{x_{r+2}}{x_1^{r+2-j}}\right)+\frac{x_{r+2}}{x_1^{r+2-i}(1-x_1)^{i-1}}
\label{Eq:ui}
\end{align}
for $i=2,\ldots,r+2$. Then we have
$$\Phi^t(1;q)=\Psi^t_n(q).$$
\end{lem}

\proof
This proof is similar as but much simpler than that of \cite[Lemma 18]{LW}. Hence we give the proof here. Using \eqref{Eq:tqFMPL-tqFMHS}, we find
\begin{align*}
\Phi^t(1)=&1+\sum\limits_{k,h_1,\ldots,h_r\geqslant  0,l>0}\sum\limits_{(k_1,\ldots,k_l)\in I(k,l,h_1,\ldots,h_r)}\sum\limits_{1\leqslant a_j\leqslant k_j\atop j=1,\ldots,l}\left(\prod\limits_{j=1}^l\binom{k_j-1}{a_j-1}\right)\\
&\times\overline{z}_n^t(a_1,\ldots,a_l;q)x_1^{k-l-\sum_{i=1}^rh_i}x_2^{l-h_1}x_3^{h_1-h_2}\cdots x_{r+1}^{h_{r-1}-h_r}x_{r+2}^{h_r}\\
=&1+\sum\limits_{k,h_1,\ldots,h_r\geqslant  0,l>0}\sum\limits_{(a_1,\ldots,a_l)\in I(k,l,h_1,\ldots,h_r)}\sum\limits_{k_j\geqslant a_j\atop j=1,\ldots,l}\left(\prod\limits_{j=1}^l\binom{k_j-1}{a_j-1}\right)\\
&\times\overline{z}_n^t(a_1,\ldots,a_l;q)x_1^{k'-l-\sum_{i=1}^rh_i'}x_2^{l-h_1'}x_3^{h_1'-h_2'}\cdots x_{r+1}^{h_{r-1}'-h_r'}x_{r+2}^{h_r'},
\end{align*}
where $(k_1,\ldots,k_l)\in I(k',l,h_1',\ldots,k_l')$ in the right-hand side of the above identity. Since
$$h_1'+\cdots+h_r'=rl-r(l-h_1')-(r-1)(h_1'-h_2')-\cdots-(h_{r-1}'-h_r')$$
and
$$h_r'=l-(l-h_1')-(h_1'-h_2')-\cdots-(h_{r-1}'-h_r'),$$
we have
\begin{align*}
&x_1^{k'-l-\sum_{i=1}^rh_i'}x_2^{l-h_1'}x_3^{h_1'-h_2'}\cdots x_{r+1}^{h_{r-1}'-h_r'}x_{r+2}^{h_r'}\\
=&x_1^{k'-(r+1)l}x_{r+2}^l\left(\frac{x_1^rx_2}{x_{r+2}}\right)^{l-h_1'}\left(\frac{x_1^{r-1}x_3}{x_{r+2}}\right)^{h_1'-h_2'}\cdots\left(\frac{x_1x_{r+1}}{x_{r+2}}\right)^{h_{r-1}'-h_r'}.
\end{align*}
Using the facts
$$l-h_1'=\sum\limits_{i=1}^l\delta_{k_i,1},\quad h_{j-1}'-h_j'=\sum\limits_{i=1}^l\delta_{k_i,j},\quad j=2,\ldots,r$$
and $k'=k_1+\cdots+k_l$, we get
\begin{align*}
&x_1^{k'-l-\sum_{i=1}^rh_i'}x_2^{l-h_1'}x_3^{h_1'-h_2'}\cdots x_{r+1}^{h_{r-1}'-h_r'}x_{r+2}^{h_r'}\\
=&\left(\frac{x_{r+2}}{x_1^{r+1}}\right)^l\prod\limits_{i=1}^lx_1^{k_i}\left(\frac{x_1^rx_2}{x_{r+2}}\right)^{\delta_{k_i,1}}\left(\frac{x_1^{r-1}x_3}{x_{r+2}}\right)^{\delta_{k_i,2}}
\cdots\left(\frac{x_1x_{r+1}}{x_{r+2}}\right)^{\delta_{k_i,r}}.
\end{align*}
Denote by
$$S_j=\sum\limits_{i\geqslant a_j}\binom{i-1}{a_j-1}x_1^{i}\left(\frac{x_1^rx_2}{x_{r+2}}\right)^{\delta_{i,1}}\left(\frac{x_1^{r-1}x_3}{x_{r+2}}\right)^{\delta_{i,2}}
\cdots\left(\frac{x_1x_{r+1}}{x_{r+2}}\right)^{\delta_{i,r}}.$$
We get
\begin{align*}
\Phi^t(1)=&1+\sum\limits_{k,h_1,\ldots,h_r\geqslant  0,l>0}\sum\limits_{(a_1,\ldots,a_l)\in I(k,l,h_1,\ldots,h_r)}\\
&\quad\times\left(\frac{x_{r+2}}{x_1^{r+1}}\right)^lS_1\cdots S_l\overline{z}_n^t(a_1,\ldots,a_l;q).
\end{align*}
For $j=1,\ldots,l$, if $a_j=m$, we have
$$S_j=\begin{cases}
X_m & \text{if\;} 1\leqslant m\leqslant r,\\
\left(\frac{x_1}{1-x_1}\right)^{m} & \text{if \;} m\geqslant r+1,
\end{cases}$$
where
$$X_m=\sum\limits_{i=m}^r\binom{i-1}{m-1}x_1^i\frac{x_1^{r+1-i}x_{i+1}-x_{r+2}}{x_{r+2}}+\left(\frac{x_1}{1-x_1}\right)^m.$$
Thus we get
\begin{align*}
\Phi^t(1)=&1+\sum\limits_{k,h_1,\ldots,h_r\geqslant  0,l>0}\sum\limits_{(a_1,\ldots,a_l)\in I(k,l,h_1,\ldots,h_r)}X_1^{l-h_1}X_2^{h_1-h_2}\cdots X_r^{h_{r-1}-h_r}\\
&\quad\times \left(\frac{x_{r+2}}{x_1^{r+1}}\right)^l\left(\frac{x_1}{1-x_1}\right)^{a_1+\cdots+a_l-[(l-h_1)+2(h_1-h_2)+\cdots+r(h_{r-1}-h_r)]}\\
&\quad\times\overline{z}_n^t(a_1,\ldots,a_l;q)\\
=&1+\sum\limits_{k,h_1,\ldots,h_r\geqslant  0,l>0}\sum\limits_{(a_1,\ldots,a_l)\in I(k,l,h_1,\ldots,h_r)}\left(\frac{x_1}{1-x_1}\right)^{k}\left(\frac{x_{r+2}}{x_1^{r+1}}\right)^l\\
&\;\times \left(X_1\frac{1-x_1}{x_1}\right)^{l-h_1}\left(X_2\left(\frac{1-x_1}{x_1}\right)^2\right)^{h_1-h_2}\cdots\left(X_r\left(\frac{1-x_1}{x_1}\right)^r\right)^{h_{r-1}-h_r}\\
&\quad\times\overline{z}_n^t(a_1,\ldots,a_l;q).
\end{align*}
As
$$k=k-l-\sum\limits_{i=1}^rh_i+(l-h_1)+2(h_1-h_2)+\cdots+r(h_{r-1}-h_r)+(r+1)h_r$$
and
$$l=(l-h_1)+(h_1-h_2)+\cdots+(h_{r-1}-h_r)+h_r,$$
we obtain
\begin{align*}
\Phi^t(1)=&1+\sum\limits_{k,h_1,\ldots,h_r\geqslant  0,l>0}\sum\limits_{(a_1,\ldots,a_l)\in I(k,l,h_1,\ldots,h_r)}\left(\frac{x_1}{1-x_1}\right)^{k-l-\sum\limits_{i=1}^rh_i}\\
&\quad\times \left(\frac{x_{r+2}}{x_1^{r+1}}X_1\right)^{l-h_1}\left(\frac{x_{r+2}}{x_1^{r+1}}X_2\right)^{h_1-h_2}\cdots\left(\frac{x_{r+2}}{x_1^{r+1}}X_r\right)^{h_{r-1}-h_r}\\
&\quad\times\left(\frac{x_{r+2}}{(1-x_1)^{r+1}}\right)^{h_r}\overline{z}_n^t(a_1,\ldots,a_l;q),
\end{align*}
from which the result follows smoothly.
\qed

Now we return to the proof of Theorem \ref{Thm:GeneratingFunction}.

\noindent {\bf Proof of Theorem \ref{Thm:GeneratingFunction}.}
It is easy to see that
$$\begin{cases}
u_1=\frac{x_1}{1-x_1}, &\\
u_{r+2}=\frac{x_{r+2}}{(1-x_1)^{r+1}} &
\end{cases}\qquad\text{and}\qquad
\begin{cases}
x_1=\frac{u_1}{1+u_1}, &\\
x_{r+2}=\frac{u_{r+2}}{(1+u_1)^{r+1}}
\end{cases}$$
are equivalent. Hence from Lemma \ref{Lem:Phi-1-Psi}, it is sufficient to show that \eqref{Eq:xi} and \eqref{Eq:ui} are equivalent for $i=2,\ldots,r+1$. In fact, for $i=2,\ldots,r+1$, \eqref{Eq:ui} are equivalent to
\begin{align}
\begin{pmatrix}
u_2\\
u_3\\
\vdots\\
u_{r+1}
\end{pmatrix}=T\begin{pmatrix}
x_2-\frac{x_{r+2}}{x_1^r}\\
x_3-\frac{x_{r+2}}{x_1^{r-1}}\\
\vdots\\
x_{r+1}-\frac{x_{r+2}}{x_1}
\end{pmatrix}+\begin{pmatrix}
\frac{x_{r+2}}{x_1^r(1-x_1)}\\
\frac{x_{r+2}}{x_1^{r-1}(1-x_1)^2}\\
\vdots\\
\frac{x_{r+2}}{x_1(1-x_1)^r}
\end{pmatrix},
\label{Eq:u-x-Matrix}
\end{align}
where $T=(t_{ij})_{r\times r}$ is an upper triangular matrix defined by
$$t_{ij}=\begin{cases}
\binom{j-1}{i-1} &\text{if\;} i\leqslant j,\\
0 & \text{if\;} i>j.
\end{cases}$$
Now $T$ is invertible, and the inverse of $T$ is $T^{-1}=(t_{ij}')_{r\times r}$ with
$$t_{ij}'=\begin{cases}
(-1)^{j-i}\binom{j-1}{i-1} &\text{if\;} i\leqslant j,\\
0 & \text{if\;} i>j.
\end{cases}$$
Therefore \eqref{Eq:u-x-Matrix} is equivalent to
\begin{align*}
\begin{pmatrix}
x_2\\
x_3\\
\vdots\\
x_{r+1}
\end{pmatrix}=&T^{-1}\begin{pmatrix}
u_2-\frac{x_{r+2}}{x_1^r(1-x_1)}\\
u_3-\frac{x_{r+2}}{x_1^{r-1}(1-x_1)^2}\\
\vdots\\
u_{r+1}-\frac{x_{r+2}}{x_1(1-x_1)^r}
\end{pmatrix}+\begin{pmatrix}
\frac{x_{r+2}}{x_1^r}\\
\frac{x_{r+2}}{x_1^{r-1}}\\
\vdots\\
\frac{x_{r+2}}{x_1}
\end{pmatrix}\\
=&T^{-1}\begin{pmatrix}
u_2-\frac{u_{r+2}}{u_1^r}\\
u_3-\frac{u_{r+2}}{u_1^{r-1}}\\
\vdots\\
u_{r+1}-\frac{u_{r+2}}{u_1}
\end{pmatrix}+\begin{pmatrix}
\frac{u_{r+2}}{u_1^r(1+u_1)}\\
\frac{u_{r+2}}{u_1^{r-1}(1+u_1)^2}\\
\vdots\\
\frac{u_{r+2}}{u_1(1+u_1)^r}
\end{pmatrix},
\end{align*}
which is equivalent to \eqref{Eq:xi} for $i=2,\ldots,r+1$.
\qed


\section{Specialization at roots of unity}\label{Sec:RootUnity}

Let $q=\zeta_n$ be a primitive $n$-th root of unity. In this section, we prove Corollary \ref{Cor:Sum-rational}, Theorem \ref{Thm:OhnoZagier-tqFMHS}, Corollary \ref{Cor:SumFormula-tqFHMS} and Corollary \ref{Cor:Eva-k-l}.


\subsection{Rationalities of the sums $G^t_n(k,l,h_1,\ldots,h_r;\zeta_n)$}\label{Subsec:Rational-Sum}

Assume that $\alpha_1^t,\ldots,\alpha_{r+1}^t$ are determined by
\begin{align*}
\begin{cases}
\alpha_1^t+\cdots+\alpha_{r+1}^t=-(x_1+tx_2), &\\
\sum\limits_{1\leqslant i_1<\cdots<i_j\leqslant r+1}\alpha_{i_1}^t\cdots\alpha_{i_j}^t=-t(x_{j+1}-x_1x_j), & j=2,\ldots,r+1.
\end{cases}
\end{align*}
Then the polynomial $P^t(T)$ defined in Theorem \ref{Thm:GeneratingFunction} can be factored as
$$P^t(T)=(T+\alpha_1^t)\cdots(T+\alpha_{r+1}^t).$$
Hence using the formula $\prod_{j=1}^{n-1}(T-\zeta_n^j)=\frac{T^n-1}{T-1}$, we get
\begin{align*}
\prod\limits_{j=1}^{n-1}P^t(1-\zeta_n^j)=&\prod\limits_{j=1}^{n-1}\prod\limits_{i=1}^{r+1}(1+\alpha_i^t-\zeta_n^j)=\prod\limits_{i=1}^{r+1}\frac{(1+\alpha_i^t)^n-1}{\alpha_i^t}\\
=&-\frac{\prod\limits_{i=1}^{r+1}\left[(1+\alpha_i^t)^n-1\right]}{t(x_{r+2}-x_1x_{r+1})}.
\end{align*}
Therefore following from Theorem \ref{Thm:GeneratingFunction}, we find
\begin{align}
\Psi^t_n(\zeta_n)=\frac{\prod\limits_{i=1}^{r+1}\left[(1+\alpha_i^{t-1})^n-1\right]}{\prod\limits_{i=1}^{r+1}\left[(1+\alpha_i^t)^n-1\right]}\frac{t}{t-1}.
\label{Eq:n-root}
\end{align}

Let $u$ be a variable. Let us consider the generating function
\begin{align*}
&\sum\limits_{n=1}^\infty \frac{\prod\limits_{i=1}^{r+1}\left[(1+\alpha_i^{t})^n-1\right]}{n}u^n\\
=&\sum\limits_{n=1}^\infty\frac{u^n}{n}\left\{(-1)^{r+1}+\sum\limits_{j=1}^{r+1}(-1)^{r+1-j}\sum\limits_{1\leqslant i_1<\cdots<i_j\leqslant r+1}\left[(1+\alpha_{i_1}^t)\cdots(1+\alpha_{i_j}^t)\right]^n\right\}\\
=&(-1)^r\log(1-u)+\sum\limits_{j=1}^{r+1}(-1)^{r-j}\sum\limits_{1\leqslant i_1<\cdots<i_j\leqslant r+1}\log\left[1-(1+\alpha_{i_1}^t)\cdots(1+\alpha_{i_j}^t)u\right].
\end{align*}
Define
\begin{align*}
&F_{r,0}^t=1-u,\\
&F_{r,j}^t=\prod\limits_{1\leqslant i_1<\cdots<i_j\leqslant r+1}\left[1-(1+\alpha_{i_1}^t)\cdots(1+\alpha_{i_j}^t)u\right],\quad j=1,\ldots,r+1,
\end{align*}
which are polynomials in $u$. We have
\begin{align}
\sum\limits_{n=1}^\infty \frac{\prod\limits_{i=1}^{r+1}\left[(1+\alpha_i^{t})^n-1\right]}{n}u^n=(-1)^r\log\left[\prod\limits_{j=0}^{r+1}\left(F_{r,j}^t\right)^{(-1)^j}\right].
\label{Eq:Generating-alpha}
\end{align}

Since the coefficients of the powers of $u$ in $F_{r,j}^t$ are symmetric polynomials in $\alpha_1^t,\ldots,\alpha_{r+1}^t$ with integer coefficients, we get
$$\prod\limits_{i=1}^{r+1}\left[(1+\alpha_i^{t})^n-1\right]\in\mathbb{Q}[t][[u_1,\ldots,u_{r+2}]].$$
And hence we have
$$G^t_n(k,l,h_1,\ldots,h_r;\zeta_n)\in\mathbb{Q}[t]$$
for any nonnegative integers $k,l,h_1,\ldots,h_r$. Then Corollary \ref{Cor:Sum-rational} is proved.

\subsection{The case of $r=1$}

Let $r=1$. Since
$$x_1=\frac{u_1}{1+u_1},\quad x_2=u_2-\frac{u_3}{1+u_1},\quad x_3=\frac{u_3}{(1+u_1)^2},$$
we get
$$\begin{cases}
\alpha_1^t+\alpha_2^t=-\frac{u_1}{1+u_1}-t\frac{u_2+u_1u_2-u_3}{1+u_1},&\\
\alpha_1^t\alpha_2^t=-t\frac{u_3-u_1u_2}{1+u_1}.&
\end{cases}$$
Thus we have
\begin{align*}
&F_{1,1}^t=1-\frac{2+u_1-t(u_2+u_1u_2-u_3)}{1+u_1}u+\frac{1-tu_2}{1+u_1}u^2,\\
&F_{1,2}^t=1-\frac{1-tu_2}{1+u_1}u.
\end{align*}
From \eqref{Eq:Generating-alpha}, we find
\begin{align*}
&\sum\limits_{n=1}^\infty\frac{\prod\limits_{i=1}^{2}\left[(1+\alpha_i^{t})^n-1\right]}{n}u^n=\log\left(1-\frac{t\frac{u_3-u_1u_2}{1+u_1}u}{(1-u)\left(1-\frac{1-tu_2}{1+u_1}u\right)}\right)\\
=&-\sum\limits_{m=1}^\infty\frac{1}{m}\frac{t^m(u_3-u_1u_2)^mu^m}{(1+u_1)^m(1-u)^m\left(1-\frac{1-tu_2}{1+u_1}u\right)^m}.
\end{align*}
Using the formula
$$(1-u)^{-m}=\sum\limits_{a\geqslant 0}\binom{m+a-1}{a}u^a,$$
we get
\begin{align*}
&\sum\limits_{n=1}^\infty\frac{\prod\limits_{i=1}^{2}\left[(1+\alpha_i^{t})^n-1\right]}{n}u^n=-\sum\limits_{m\geqslant 1,a,b\geqslant 0}\frac{1}{m}\binom{m+a-1}{a}\binom{m+b-1}{b}\\
&\qquad\times t^m(u_3-u_1u_2)^m(1-tu_2)^b(1+u_1)^{-m-b}u^{a+b+m},
\end{align*}
which implies that
\begin{align*}
\prod\limits_{i=1}^{2}\left[(1+\alpha_i^{t})^n-1\right]=&-\sum\limits_{a,b\geqslant 0\atop a+b\leqslant n-1}\frac{n}{n-a-b}\binom{n-b-1}{a}\binom{n-a-1}{b}\\
&\quad\times t^{n-a-b}(u_3-u_1u_2)^{n-a-b}(1-tu_2)^b(1+u_1)^{a-n}.
\end{align*}
Then Theorem \ref{Thm:OhnoZagier-tqFMHS} follows from \eqref{Eq:n-root} and the above formula.


\subsection{The sum formula}

Taking $u_3=u_1u_2$ in Theorem \ref{Thm:OhnoZagier-tqFMHS}, we have
\begin{align*}
U_n^t(u_1,u_2,u_1u_2)=&\sum\limits_{a,b\geqslant 0\atop a+b=n-1}(1+u_1)^a(1-tu_2)^b\\
=&\sum_{b=0}^{n-1}\sum\limits_{i=0}^b\sum\limits_{j=0}^{n-1-b}\binom{b}{i}\binom{n-1-b}{j}u_1^j(-tu_2)^i.
\end{align*}
Changing the order of the summations, we get
$$U_n^t(u_1,u_2,u_1u_2)=\sum\limits_{0\leqslant i,j\leqslant n-1}\left(\sum\limits_{b=i}^{n-1-j}\binom{b}{i}\binom{n-1-b}{j}\right)u_1^j(-tu_2)^i.$$
Using the Chu-Vandermonde identity
$$\sum\limits_{m=0}^n\binom{m}{j}\binom{n-m}{k-j}=\binom{n+1}{k+1},\quad (1\leqslant j\leqslant k\leqslant n),$$
we find
$$U_n^t(u_1,u_2,u_1u_2)=\sum\limits_{0\leqslant i,j\leqslant n-1}\binom{n}{i+j+1}u_1^j(-tu_2)^i.$$
Therefore we have
\begin{align}
\Psi_n^t(u_1,u_2,u_1u_2;\zeta_n)=\frac{\sum\limits_{0\leqslant i,j\leqslant n-1}\binom{n}{i+j+1}u_1^j((1-t)u_2)^i}{\sum\limits_{0\leqslant i,j\leqslant n-1}\binom{n}{i+j+1}u_1^j(-tu_2)^i}.
\label{Eq:Sum-root}
\end{align}

On the one hand, since
$$U_n^t(u_1,u_2,u_1u_2)=n+\sum\limits_{0\leqslant i,j\leqslant n-1}\binom{n}{i+j+1}(1-\delta_{i,0}\delta_{j,0})u_1^j(-tu_2)^i,$$
from \eqref{Eq:Sum-root} we find
\begin{align*}
&\Psi_n^t(u_1,u_2,u_1u_2;\zeta_n)=\frac{\frac{1}{n}\sum\limits_{0\leqslant i,j\leqslant n-1}\binom{n}{i+j+1}u_1^j((1-t)u_2)^i}{1+\frac{1}{n}\sum\limits_{0\leqslant i,j\leqslant n-1}\binom{n}{i+j+1}(1-\delta_{i,0}\delta_{j,0})u_1^j(-tu_2)^i}\\
=&\sum\limits_{m=0}^\infty\frac{(-1)^m}{n^{m+1}}\sum\limits_{0\leqslant i_a,j_a\leqslant n-1\atop a=0,\ldots,m}\prod\limits_{a=0}^m\binom{n}{i_a+j_a+1}\prod\limits_{a=1}^m(1-\delta_{i_a,0}\delta_{j_a,0})\\
&\quad\times(1-t)^{i_0}(-t)^{i_1+\cdots+i_m}u_1^{j_0+\cdots+j_m}u_2^{i_0+\cdots+i_m}\\
=&\sum\limits_{m=0}^\infty\frac{(-1)^m}{n^{m+1}}\sum\limits_{0\leqslant i_a\leqslant j_a\leqslant n-1\atop a=0,\ldots,m}\prod\limits_{a=0}^m\binom{n}{j_a+1}\prod\limits_{a=1}^m(1-\delta_{i_a,0}\delta_{j_a,0})\\
&\quad\times(1-t)^{i_0}(-t)^{i_1+\cdots+i_m}u_1^{j_0+\cdots+j_m-i_0-\cdots-i_m}u_2^{i_0+\cdots+i_m}.
\end{align*}
By the definition, we have
$$\Psi_n^t(u_1,u_2,u_1u_2;\zeta_n)=\sum\limits_{k,l\geqslant 0}G_n^t(k,l;\zeta_n)u_1^{k-l}u_2^l.$$
Comparing the coefficients of $u_1^{k-l}u_2^l$ for $k\geqslant l\geqslant 0$, we get
\begin{align*}
G_n^t(k,l;\zeta_n)=&\sum\limits_{m=0}^\infty\frac{(-1)^m}{n^{m+1}}\sum\limits_{{i_0+\cdots+i_m=l\atop j_0+\cdots+j_m=k}\atop {0\leqslant i_a\leqslant j_a\leqslant n-1\atop a=0,\ldots,m}}\prod\limits_{a=0}^m\binom{n}{j_a+1}\\
&\quad\times\prod\limits_{a=1}^m(1-\delta_{i_a,0}\delta_{j_a,0})(1-t)^{i_0}(-t)^{l-i_0}\\
=&\sum\limits_{m=0}^k\frac{(-1)^m}{n^{m+1}}\sum\limits_{{{i_0+\cdots+i_m=l\atop j_0+\cdots+j_m=k}\atop {0\leqslant i_a\leqslant j_a\leqslant n-1\atop a=0,\ldots,m}}\atop j_1,\ldots,j_m\geqslant 1}\prod\limits_{a=0}^m\binom{n}{j_a+1}(1-t)^{i_0}(-t)^{l-i_0},
\end{align*}
which is just \eqref{Eq:SumFormula-1}.

On the other hand, we have
$$U_n^t(u_1,u_2,u_1u_2)=\sum\limits_{j=0}^{n-1}\binom{n}{j+1}u_1^j+\sum\limits_{i=1}^{n-1}\left(\sum\limits_{j=i}^{n-1}\binom{n}{j+1}u_1^{j-i}\right)(-tu_2)^i,$$
where
$$\sum\limits_{j=0}^{n-1}\binom{n}{j+1}u_1^j=\frac{(1+u_1)^n-1}{u_1}.$$
For $i=0,\ldots,n-1$, set
$$f_i(u_1)=\frac{u_1}{(1+u_1)^n-1}\sum\limits_{j=i}^{n-1}\binom{n}{j+1}u_1^{j-i}.$$
Then from \eqref{Eq:Sum-root} we get
\begin{align*}
&\Psi_n^t(u_1,u_2,u_1u_2;\zeta_n)=\frac{\sum\limits_{i=0}^{n-1}f_i(u_1)((1-t)u_2)^i}{1+\sum\limits_{i=1}^{n-1}f_i(u_1)(-tu_2)^i}\\
=&\sum\limits_{{m\geqslant 0\atop 0\leqslant i_0\leqslant n-1}\atop{1\leqslant i_a\leqslant n-1\atop a=1,\ldots,m}}(-1)^mf_{i_0}(u_1)\cdots f_{i_m}(u_1)(1-t)^{i_0}(-t)^{i_1+\cdots+i_m}u_2^{i_0+\cdots+i_m}.
\end{align*}
By \cite[(3.13)]{BTT18},
$$\frac{u_1}{(1+u_1)^n-1}=\frac{1}{n}\left(1-\sum\limits_{l=1}^\infty\overline{z}_n(l;\zeta_n)u_1^l\right)=-\frac{1}{n}\sum\limits_{l=0}^\infty\overline{z}_n(l;\zeta_n)u_1^l,$$
with the convention $\overline{z}_n(0;\zeta_n)=-1$. Therefore, we have
$$f_i(u_1)=-\frac{1}{n}\sum\limits_{l\geqslant 0,i\leqslant j\leqslant n-1}\binom{n}{j+1}\overline{z}_n(l;\zeta_n)u_1^{l+j-i},$$
which gives
\begin{align*}
&\Psi_n^t(u_1,u_2,u_1u_2;\zeta_n)=-\sum\limits_{m=0}^\infty\frac{1}{n^{m+1}}\sum\limits_{{0\leqslant i_0\leqslant j_0\leqslant n-1\atop 1\leqslant i_a\leqslant j_a\leqslant n-1}\atop {a=1,\ldots,m\atop l_0,\ldots,l_m\geqslant 0}}\prod\limits_{a=0}^m\binom{n}{j_a+1}\overline{z}_n(l_a;\zeta_n)\\
&\qquad\times (1-t)^{i_0}(-t)^{i_1+\cdots+i_m}u_1^{l_0+\cdots+l_m+j_0+\cdots+j_m-i_0-\cdots-i_m}u_2^{i_0+\cdots+i_m}.
\end{align*}
Comparing the coefficients, we prove \eqref{Eq:SumFormula-2}.

\subsection{The values $\overline{z}_n^t(\{k\}^l;\zeta_n)$}

Setting $u_1=u_3=0$ in Theorem \ref{Thm:OhnoZagier-tqFMHS}, we get
$$\sum\limits_{l=0}^\infty\overline{z}_n^t(\{1\}^l;\zeta_n)u_2^l=\frac{\sum\limits_{a,b\geqslant 0\atop a+b=n-1}(1-(t-1)u_2)^b}{\sum\limits_{a,b\geqslant 0\atop a+b=n-1}(1-tu_2)^b}.$$
Now since
\begin{align*}
\sum\limits_{a,b\geqslant 0\atop a+b=n-1}(1-tu_2)^b=&\frac{(1-tu_2)^n-1}{-tu_2}=\sum\limits_{i=0}^{n-1}\binom{n}{i+1}(-tu_2)^i\\
=&n+\sum\limits_{i=1}^{n-1}\binom{n}{i+1}(-tu_2)^i,
\end{align*}
we have
\begin{align*}
&\sum\limits_{l=0}^\infty\overline{z}_n^t(\{1\}^l;\zeta_n)u_2^l=\frac{\frac{1}{n}\sum\limits_{i=0}^{n-1}\binom{n}{i+1}((1-t)u_2)^i}{1+\frac{1}{n}\sum\limits_{i=1}^{n-1}\binom{n}{i+1}(-tu_2)^i}\\
=&\sum\limits_{m=0}^\infty\frac{(-1)^m}{n^{m+1}}\sum\limits_{0\leqslant i_0\leqslant n-1\atop 1\leqslant i_1,\ldots,i_m\leqslant n-1}\prod\limits_{a=0}^m\binom{n}{i_a+1}(1-t)^{i_0}(-t)^{i_1+\cdots+i_m}u_2^{i_0+\cdots+i_m}.
\end{align*}
Comparing the coefficients of $u_2^l$, we get the evaluation formula of $\overline{z}_n^t(\{1\}^l;\zeta_n)$ as displayed in \eqref{Eq:Eva-1-l}.

Setting $u_1=u_2=0$ in Theorem \ref{Thm:OhnoZagier-tqFMHS}, we find
$$\sum\limits_{l=0}^\infty\overline{z}_n^t(\{2\}^l;\zeta_n)u_3^l=\frac{U_n^{t-1}(0,0,u_3)}{U_n^t(0,0,u_3)},$$
where
\begin{align*}
U_n^t(0,0,u_3)=&\sum\limits_{a,b\geqslant 0\atop a+b\leqslant n-1}\frac{1}{n-a-b}\binom{n-a-1}{b}\binom{n-b-1}{a}(tu_3)^{n-a-b-1}\\
=&\sum\limits_{i=0}^{n-1}\left(\sum\limits_{a,b\geqslant 0\atop a+b=i}\binom{n-a-1}{b}\binom{n-b-1}{a}\right)\frac{1}{n-i}(tu_3)^{n-i-1}.
\end{align*}
Now since
\begin{align*}
&(1-T)^{-(n-i)}=\sum\limits_{a=0}^\infty\binom{n-i-1+a}{a}T^a,\\
&(1-T)^{-2(n-i)}=\sum\limits_{a=0}^\infty\binom{2n-2i-1+a}{a}T^a,
\end{align*}
we find
$$\sum\limits_{a,b\geqslant 0\atop a+b=i}\binom{n-a-1}{b}\binom{n-b-1}{a}=\binom{2n-i-1}{i}.$$
Therefore, we get
$$U_n^t(0,0,u_3)=\sum\limits_{i=0}^{n-1}\frac{1}{i+1}\binom{n+i}{2i+1}(tu_3)^{i}=n+\sum\limits_{i=1}^{n-1}\frac{1}{i+1}\binom{n+i}{2i+1}(tu_3)^{i},$$
which implies that
\begin{align*}
&\sum\limits_{l=0}^\infty\overline{z}_n^t(\{2\}^l;\zeta_n)u_3^l=\frac{\frac{1}{n}\sum\limits_{i=0}^{n-1}\frac{1}{i+1}\binom{n+i}{2i+1}((t-1)u_3)^{i}}{1+\frac{1}{n}\sum\limits_{i=1}^{n-1}\frac{1}{i+1}\binom{n+i}{2i+1}(tu_3)^{i}}\\
=&\sum\limits_{m=0}^\infty\frac{(-1)^m}{n^{m+1}}\sum\limits_{0\leqslant i_0\leqslant n-1\atop 1\leqslant i_1,\ldots,i_m\leqslant n-1}\prod\limits_{a=0}^m\frac{1}{i_a+1}\binom{n+i_a}{2i_a+1}(t-1)^{i_0}t^{i_1+\cdots+i_m}u_3^{i_0+\cdots+i_m}.
\end{align*}
Comparing the coefficients of $u_3^l$, we get the evaluation formula of $\overline{z}_n^t(\{2\}^l;\zeta_n)$ as displayed in \eqref{Eq:Eva-2-l}.

More generally, setting $u_1=\cdots=u_{r+1}=0$ and $q=\zeta_n$ in Theorem \ref{Thm:GeneratingFunction}, we get
$$\sum\limits_{l=0}^\infty \overline{z}_n^t(\{r+1\}^l;\zeta_n)u_{r+2}^l=\frac{\prod\limits_{j=1}^{n-1}P^{t-1}(1-\zeta_n^j)}{\prod\limits_{j=1}^{n-1}P^t(1-\zeta_n^j)}.$$
To compute the polynomial $P^t(T)$ in this case, we need the following lemma.

\begin{lem}\label{Lem:u1-ur-0-xi}
If $u_1=\cdots=u_{r+1}=0$, then for any integer $i$ with $1\leqslant i\leqslant r+2$, we have
$$x_i=(-1)^{r-i}\binom{r}{r+2-i}u_{r+2}.$$
\end{lem}

\proof
By \eqref{Eq:x1}, we have $x_1=0$. By \eqref{Eq:xi}, we have $x_{r+2}=u_{r+2}$ and
\begin{align*}
\frac{x_i}{u_{r+2}}=&\lim\limits_{u_1\rightarrow 0}\frac{1-(1+u_1)^{i-1}\sum\limits_{j=i}^{r+1}(-1)^{j-i}\binom{j-2}{i-2}u_1^{j-i}}{u_1^{r+2-i}}\\
=&\lim\limits_{u_1\rightarrow 0}\frac{1-(1+u_1)^{i-1}\sum\limits_{j=0}^{r+1-i}(-1)^{j}\binom{i+j-2}{j}u_1^{j}}{u_1^{r+2-i}}
\end{align*}
for $2\leqslant i\leqslant r+1$. Fix an integer $i$ such that $2\leqslant i\leqslant r+1$. Using the L'H\^{o}pital's rule, we find
$$\frac{x_i}{u_{r+2}}=\lim\limits_{u_1\rightarrow 0}\frac{-(1+u_1)\sum\limits_{j=1}^{r+1-i}(-1)^jj\binom{i+j-2}{j}u_1^{j-1}-(i-1)\sum\limits_{j=0}^{r+1-i}(-1)^j\binom{i+j-2}{j}u_1^j}{(r+2-i)u_1^{r+1-i}}.$$
In particular, we have
$$x_{r+1}=-ru_{r+2}.$$
Now assume that $2\leqslant i\leqslant r$. Using the L'H\^{o}pital's rule again, we get
$$\frac{x_i}{u_{r+2}}=\lim\limits_{u_1\rightarrow 0}\frac{f(u_1)}{(r+2-i)(r+1-i)u_1^{r-i}},$$
where
\begin{align*}
f(u_1)=&-\sum\limits_{j=2}^{r+1-i}(-1)^jj(j-1)\binom{i+j-2}{j}u_1^{j-2}-\sum\limits_{j=1}^{r+1-i}(-1)^jj^2\binom{i+j-2}{j}u_1^{j-1}\\
&\quad -(i-1)\sum\limits_{j=1}^{r+1-i}(-1)^jj\binom{i+j-2}{j}u_1^{j-1}.
\end{align*}
Taking the sum of the latter two summations, we find
\begin{align*}
f(u_1)=&\sum\limits_{j=1}^{r-i}(-1)^j\frac{(i+j-1)!}{(j-1)!(i-2)!}u_1^{j-1}-\sum\limits_{j=1}^{r+1-i}(-1)^j\frac{(i+j-1)!}{(j-1)!(i-2)!}u_1^{j-1}\\
=&(-1)^{r-i}\frac{r!}{(r-i)!(i-2)!}u_1^{r-i},
\end{align*}
which implies the desired formula.
\qed

Using Lemma \ref{Lem:u1-ur-0-xi}, if $u_1=\cdots=u_{r+1}=0$, we have
$$P^t(T)=T^{r+1}-t(1-T)^ru_{r+2}.$$
Hence from Theorem \ref{Thm:GeneratingFunction}, we get the following corollary. Here the result for $k=1$ is followed from setting $r=1$ and $u_1=u_3=0$ in Theorem \ref{Thm:GeneratingFunction}.

\begin{cor}
Let $k\in\mathbb{N}$ and $v$ be a variable. Then we have
$$\sum\limits_{l=0}^\infty\overline{z}_n^t(\{k\}^l;\zeta_n)v^l=\frac{\prod\limits_{j=1}^{n-1}\widetilde{P}^{t-1}(\zeta_n^j)}{\prod\limits_{j=1}^{n-1}\widetilde{P}^{t}(\zeta_n^j)},$$
where
$$\widetilde{P}^{t}(T)=\begin{cases}
(1-T)^2-t(1-T)v & \text{if\;} k=1,\\
(1-T)^{k}-tT^{k-1}v & \text{if\;} k\geqslant 2.
\end{cases}$$
\end{cor}

Now assume that $k\geqslant 2$. Set
$$\widetilde{P}^{t}(T)=(\beta_1^t-T)\cdots(\beta_k^t-T).$$
Then the elementary symmetric polynomials of $\beta_1^t,\ldots,\beta_k^t$ are given by
\begin{align}
e_j=\sum\limits_{1\leqslant i_1<\cdots<i_j\leqslant k}\beta_{i_1}^t\cdots\beta_{i_j}^t=\begin{cases}
k+(-1)^ktv & \text{if\;} j=1,\\
\binom{k}{j} & \text{if\;} j=2,\ldots,k.
\end{cases}
\label{Eq:Beta}
\end{align}
Similar as the discussion in Subsection \ref{Subsec:Rational-Sum}, we get the following result.

\begin{cor}\label{Cor:Generating-k-l}
Let $k\geqslant 2$ be an integer and $u,v$ be variables. Let $\beta_1^t,\ldots,\beta_k^t$ be determined by \eqref{Eq:Beta}.
Let $\widetilde{F}_0^t=1-u$ and
$$\widetilde{F}_j^t=\prod\limits_{1\leqslant i_1<\cdots<i_j\leqslant k}(1-\beta_{i_1}^t\cdots\beta_{i_j}^tu),\quad (1\leqslant j\leqslant k).$$
Then we have
$$\sum\limits_{l=0}^\infty\overline{z}_n^t(\{k\}^l;\zeta_n)v^l=\frac{H_n^{t-1}}{H_n^t}\frac{t}{t-1},$$
where $H_n^t$ is determined by
$$\sum\limits_{n=1}^\infty\frac{H_n^t}{n}u^n=(-1)^{k-1}\log\left[\prod\limits_{j=0}^{k}\left(\widetilde{F}_j^t\right)^{(-1)^j}\right].$$
\end{cor}

Applying Corollary \ref{Cor:Generating-k-l} and with some computations, one can obtain the evaluation formulas of $\overline{z}_n^t(\{k\}^l;\zeta_n)$ for any $k\geqslant 2$. Here we take $k=3$ as an example to prove \eqref{Eq:Eva-3-l}.

In this case, we have
$$e_1=3-tv,\quad e_2=3,\quad e_3=1.$$
Hence we find
$$\widetilde{F}_1^t=(1-u)^3+tvu,\quad \widetilde{F}_2^t=(1-u)^3-tvu^2,\quad \widetilde{F}_3^t=1-u,$$
which implies
\begin{align*}
\sum\limits_{n=1}^\infty\frac{H_n^t}{n}u^n=&\log\frac{(1-u)^3-tvu^2}{(1-u)^3+tvu}=\log\left(1-\frac{tvu^2}{(1-u)^3}\right)-\log\left(1+\frac{tvu}{(1-u)^3}\right)\\
=&\sum\limits_{i=1}^\infty\frac{-(tvu^2)^i+(-tvu)^i}{i(1-u)^{3i}}\\
=&\sum\limits_{i=1}^\infty\sum\limits_{a=0}^\infty\frac{-(tvu^2)^i+(-tvu)^i}{i}\binom{3i+a-1}{a}u^a.
\end{align*}
Therefore, we have
$$\frac{H_n^t}{n}=-\sum\limits_{i=0}^{n-1}\frac{1}{i+1}\left[\binom{n+i}{3i+2}+(-1)^i\binom{n+2i+1}{3i+2}\right](tv)^{i+1},$$
which implies that
\begin{align*}
\sum\limits_{l=0}^\infty\overline{z}_n^t(\{3\}^l;\zeta_n)v^l=&\frac{\sum\limits_{i=0}^{n-1}\frac{1}{i+1}\left[\binom{n+i}{3i+2}+(-1)^i\binom{n+2i+1}{3i+2}\right]((t-1)v)^{i}}
{\sum\limits_{i=0}^{n-1}\frac{1}{i+1}\left[\binom{n+i}{3i+2}+(-1)^i\binom{n+2i+1}{3i+2}\right](tv)^{i}}\\
=&\frac{\frac{1}{n^2}\sum\limits_{i=0}^{n-1}\frac{1}{i+1}\left[\binom{n+i}{3i+2}+(-1)^i\binom{n+2i+1}{3i+2}\right]((t-1)v)^{i}}
{1+\frac{1}{n^2}\sum\limits_{i=1}^{n-1}\frac{1}{i+1}\left[\binom{n+i}{3i+2}+(-1)^i\binom{n+2i+1}{3i+2}\right](tv)^{i}}.
\end{align*}
We then prove \eqref{Eq:Eva-3-l} without difficulty.



\end{document}